\documentclass{article}
\usepackage{indentfirst,latexsym,bm,amsmath,amsthm,amsfonts}
\usepackage[dvips]{graphicx}
\def\reg{\textrm{reg}}
\begin{document}
\title{Multiplicities of representations in spaces \\of modular forms}

\author{BinYong Hsie
 \\LMAM, Department of Mathematics, PeKing University,\\BeiJing, P.R.China
      100871
\\ e-mail: byhsie@math.pku.edu.cn}

\date{}
\maketitle{}

\begin{abstract}This paper shows that for a given irreducible
representation $\rho$ of $\Gamma/\Gamma_1$, the two functions
dim($M_k(\Gamma_1,\rho)$) and dim($S_k(\Gamma_1,\rho)$) of $k$ are
almost linear functions.
\end{abstract}
\footnotetext
%"*"
{Mathematical Subject Classification: 11F11. Key Words: modular
form, cusp form.}

\vskip 10pt

Let $\Gamma$ and $\Gamma_1$ be two Fuchsian subgroups of
SL$_2(\mathbb{Q})$ of the first kind, with $\Gamma_1$ a normal
subgroup of $\Gamma$ of index $\mu$. Let $\Bbb H^*$ denote the
extended upper half-plane \cite{Sh} and put $Y=\Gamma\backslash
\Bbb H^*$ and $X=\Gamma_1\backslash \Bbb H^*$. Let $A_k(\ast)$
denotes the space of meromorphic modular forms for $\ast$ of
weight $k$. We define a representation $\pi_k$ of $\Gamma$ on
$A_k(\Gamma_1)$ by the formula
$$
\pi_k(\gamma\Gamma_1)( f)=f|[\gamma^{-1}]_k \phantom{............}
(\gamma\in\Gamma , \phantom{...} f\in A_k(\Gamma_1)),
$$
where the notation $[\ast]_k$ is as in Shimura \cite{Sh}. The
representation $\pi_k$ factors through $\Gamma/\Gamma_1$ to give a
representation -- also denoted $\pi_k$ -- of the latter group. The
space $M_k(\Gamma_1)$ of holomorphic modular forms of weight $k$
for $\Gamma_1$ and the subspace $S_k(\Gamma_1)$ of cusp forms of
weight $k$ for $\Gamma_1$ are both stable under $\pi_k$, and the
resulting representations of $\Gamma/\Gamma_1$ on $M_k(\Gamma_1)$
and $S_k(\Gamma_1)$ will be denoted $\rho_k$ and $\sigma_k$
respectively.

Henceforth $\rho$ denotes an irreducible complex representation of
$\Gamma/\Gamma_1$. If $-I\in \Gamma$ then the value of $\rho$ on
the coset of $-I$ in $\Gamma/\Gamma_1$ is a scalar by Schur's
lemma, and we say that $\rho$ is {\it even} or {\it odd} according
as the scalar is 1 or $-1$. If $-I\in \Gamma_1$ then $\rho$ is
automatically even.

If $\pi$ is any finite-dimensional complex representation of
$\Gamma/\Gamma_1$ then we write $\langle\rho,\pi\rangle$ for the
multiplicity of $\rho$ in $\pi$. For example $\langle \rho,
\rho_{reg} \rangle=\dim\rho$, where $\rho_{reg}$ is the regular
representation of $\Gamma/\Gamma_1$.

\vskip 3pt

{\bf Theorem} {\it Fix an irreducible complex representation
$\rho$ of $\Gamma/\Gamma_1$, and put
$$
c=\dfrac{1}{4\pi}\int_{\Gamma\backslash \Bbb H} \dfrac{dxdy}{y^2}.
$$
If $-I\notin \Gamma$ then
$$
\lim_{k\to \infty} \dfrac{\langle \rho,
\rho_k\rangle}{k\langle\rho,\rho_{reg}\rangle}= \lim_{k\to \infty}
\dfrac{\langle \rho,
\sigma_k\rangle}{k\langle\rho,\rho_{reg}\rangle}= c.
$$
If $-I\in \Gamma$ then the same assertion holds provided $k$ runs
through positive integers of the same parity as $\rho$.}

\vskip 3pt

{\bf Proof.} We shall prove our assertion only for $\rho_{k}$. The
proof for $\sigma_{k}$ is similar.

For any positive integer $k$, let $i(k)=(-1)^{k}$. If $-I\in
\Gamma_{1}$, we assume $k$ is always even.

Let $\mathcal S$ be the set of the irreducible representations of
$\Gamma/\Gamma_{1}$ with the same parity as $k$ if $-I \in
\Gamma$, and the set of all irreducible representations of
$\Gamma/\Gamma_{1}$ otherwise.

Let $p$ be a non-cusp point of $Y$ which has $\mu$ different
liftings $p_{1},...,p_{\mu}$ on $X$, and let $\tilde p_{j}$ be a
lifting of $p_{j}$ in $\Bbb H^{*}$($1\le j\le\mu$). By
Riemann-Roch theorem we know that there exists an $f\in
A_{i(k)}(\Gamma_{1})$ such that $f$ has poles only at cusps and
$f(\tilde p_{1})=1$, $f(\tilde p_{j})=0$ for all $j=2,\cdots,\mu$.
For any $\alpha\in\Gamma/\Gamma_{1}$, put
$f_{\alpha}=f|[\alpha]_{i(k)}$. For a fixed sufficiently large
even integer $n_{0}$, we may choose $\varphi=\varphi_{n_0}\in
S_{n_{0}}(\Gamma)$ such that $\varphi
f_{\alpha}\;(\alpha\in\Gamma/\Gamma_{1})$ are all cusp forms. Let
$W$ denote the $\Bbb C$-linear subspace of
$M_{i(k)+n_{0}}(\Gamma_{1})$ spanned by $\{\varphi
f_{\alpha}|\alpha\in\Gamma/\Gamma_{1}\}$. Clearly, $W$ is stable
under $\pi_{i(k)+n_{0}}$. It is easy to see that
$$W\cong
\bigoplus\limits_{\rho\in\mathcal S}\langle\rho,\rho_{\reg}\rangle
\rho
$$
since the two hand sides have the same trace.

When $k\geq i(k)+n_{0}$, the choice of $f$ ensures that the
natural map
$$W\bigotimes_{\Bbb C}M_{k-i(k)-n_{0}}(\Gamma)\rightarrow M_{k}(\Gamma_{1})
$$
is injective. Hence
$$
\dfrac{\langle\rho,\rho_{k}\rangle}{k\langle\rho,\rho_{\reg}\rangle}
\geq\dfrac{\dim(M_{k-i(k)-n_{0}}(\Gamma))}{k}. \eqno(1)
$$
For $\rho$ and $k$ as in the Theorem, in [4] Shimura showed that
$$
\lim\limits_{k\to\infty} \frac{\dim(M_{k-i}(\Gamma))}{k}=c\eqno(2)
$$
and
$$
\lim\limits_{k\to\infty}
\frac{\dim(M_{k}(\Gamma_{1}))}{k}=c\mu.\eqno(3)
$$
From (1) and (2) we easily deduce that
$$
\liminf_{k\to\infty} \dfrac{\langle\rho,\rho_{k}\rangle}{k}\geq
c\langle\rho,\rho_{\reg}\rangle.\eqno(4)
$$
It is obvious that (3) is equivalent to
$$\lim\limits_{k\to\infty}\sum\limits_{\rho\in\mathcal S}\dim(\rho)\dfrac{\langle\rho,\rho_{k}\rangle}{k}=
\sum\limits_{\rho\in\mathcal
S}\dim(\rho)\;c\langle\rho,\rho_{\reg}\rangle.
$$
Combining this equality with (4), we get
$$\lim\limits_{k\to\infty}
\dfrac{\langle\rho,\rho_{k}\rangle}{k}=
c\langle\rho,\rho_{\reg}\rangle,
$$
as desired.\qed
%0503215.ps.gz (30kb)

{\bf Acknowledgement}. I would like to thank Professor David E.
Rohrlich for his kindly help in revising this paper.

\end{document}